\documentclass{IEEEtran4PSCC}
\usepackage{cite}
\ifCLASSINFOpdf
   \usepackage[pdftex]{graphicx}
\else
   \usepackage[dvips]{graphicx}
\fi
\usepackage[cmex10]{amsmath}
\makeatletter
\let\old@ps@headings\ps@headings
\let\old@ps@IEEEtitlepagestyle\ps@IEEEtitlepagestyle
\def\psccfooter#1{%
    \def\ps@headings{%
        \old@ps@headings%
        \def\@oddfoot{\strut\hfill#1\hfill\strut}%
        \def\@evenfoot{\strut\hfill#1\hfill\strut}%
    }%
    \def\ps@IEEEtitlepagestyle{%
        \old@ps@IEEEtitlepagestyle%
        \def\@oddfoot{\strut\hfill#1\hfill\strut}%
        \def\@evenfoot{\strut\hfill#1\hfill\strut}%
    }%
    \ps@headings%
}
\makeatother

\usepackage{booktabs} 
\usepackage{multirow}
\usepackage{rotating}
\usepackage{tabularx}
\usepackage{amsmath,amsthm,amsfonts,amssymb}
\usepackage{bm}
\usepackage{xcolor}

\theoremstyle{definition}
\newtheorem{theorem}{\textbf{Theorem}}

\newtheorem{proposition}[theorem]{\textbf{Proposition}}

\newcommand{\up}[1]{^\mathrm{#1}}
\begin{document}
%
% paper title
% Titles are generally capitalized except for words such as a, an, and, as,
% at, but, by, for, in, nor, of, on, or, the, to and up, which are usually
% not capitalized unless they are the first or last word of the title.
% Linebreaks \\ can be used within to get better formatting as desired.
% Do not put math or special symbols in the title.
\title{Energy Storage Arbitrage in Two-settlement Markets: A Transformer-Based Approach}

%% To specify the authors when (number of affiliations <= 2)
\author{
\IEEEauthorblockN{Saud Alghumayjan\IEEEauthorrefmark{1}, Jiajun Han\IEEEauthorrefmark{2}, Ningkun Zheng\IEEEauthorrefmark{1}, Ming Yi\IEEEauthorrefmark{3}, Bolun Xu\IEEEauthorrefmark{1}}
\IEEEauthorblockA{\IEEEauthorrefmark{1} Earth and Environmental Engineering, \IEEEauthorrefmark{2} Electrical Engineering, \IEEEauthorrefmark{3} Data Science Institute \\
Columbia University\\
New York, NY 10027, United States\\
{\{saa2244, jh4316, nz2343, my2826, bx2177\}@columbia.edu}}
% \and
% \IEEEauthorblockN{Author n.1 Name per Affiliation B\\ Author n.2 Name per Affiliation B}
% \IEEEauthorblockA{(Affiliation B) Department Name of Organization \\
% Name of the organization, acronyms acceptable\\
% City, Country\\
% \{email author n.1, email author n.2\}@domain (if desired)}
}

% make the title area
\maketitle

\begin{abstract}
This paper presents an integrated model for bidding energy storage in day-ahead and real-time markets to maximize profits. We show that in integrated two-stage bidding, the real-time bids are independent of day-ahead settlements, while the day-ahead bids should be based on predicted real-time prices. We utilize a transformer-based model for real-time price prediction, which captures complex dynamical patterns of real-time prices, and use the result for day-ahead bidding design. For real-time bidding, we utilize a long short-term memory-dynamic programming hybrid real-time bidding model. We train and test our model with historical data from New York State, and our results showed that the integrated system achieved promising results of almost a 20\% increase in profit compared to only bidding in real-time markets, and at the same time reducing the risk in terms of the number of days with negative profits.
\end{abstract}

\begin{IEEEkeywords}
Electricity markets, Energy storage, Machine learning
\end{IEEEkeywords}

% Use this to place sponsorships
\thanksto{\noindent The authors would like to acknowledge supports from Columbia Data Science Institute, Red River Clean Energy Inc., and King Abdulaziz City for Science and Technology (KACST).}

\section{Introduction}
\subsection{Motivation}
Energy storage participants are increasingly pivotal in electricity markets. In 2018, the Federal Energy Regulatory Commission (FERC) issued Order 841, mandating the inclusion of storage in all electricity markets~\cite{ferc}. Consequently, all system operators across the United States now permit storage to place both charge and discharge bids. Such a bidding system allows storage entities to craft their offers based on market prices. As storage capacity grows, the emphasis has transitioned from ancillary service markets, characterized by limited market capacity~\cite{sioshansi}, to wholesale energy markets~\cite{Williams, us2021form}.

United States' electricity markets follow a two-stage settlement design: a day-ahead market (DAM) and a real-time market (RTM). The DAM, covering the subsequent day's 24 trading intervals, requires bids to be placed by a specified closing time for each hour. The RTM is responsible for balancing the supply and demand deviations from day-ahead schedules. Existing literature offers multiple models for energy storage arbitrage but mostly focuses on either DAM or RTM. The most common method is model predictive control (MPC)~\cite{arnold}, but its effectiveness diminishes for long look-ahead windows, such as 24 hours~\cite{abdulla2016optimal}. Stochastic dynamic programming (SDP)~\cite{Zheng} and reinforcement learning (RL)~\cite{wang} are emerging approaches but often encounter computation difficulties. 

RTM is significantly more attractive for storage participants to arbitrage as the prices are more volatile than DAM~\cite{Milstein}. On the other hand, recent studies~\cite{qin2023role} have shown while RTM provides higher profits, storage participating in both DAM and RTM provides the best social welfare outcome. Hence, system operators like the California Independent System Operator (CAISO) have mandated storage must bid both into day-ahead and real-time markets~\cite{Byrne}. However, the newest bidding data from CAISO~\cite{Caiso} shows most participants submit unreasonably high bids day-ahead to purposefully avoid day-ahead settlements, like due to lacking good joint bidding solutions.

\subsection{Contribution}
While it is lucrative for energy storage entities to engage in both DAM and RTM, the inherent challenges stem from multi-stage decision-making amidst volatile and unpredictable electricity prices. Notwithstanding the existing research~\cite{rad2015}, our study proposes a novel and practical approach integrating machine learning to bid energy storage into both day-ahead and real-time markets. Our salient contributions are:
\begin{itemize}
    \item We propose a novel energy storage arbitrage in two-settlement markets framework that combines a transformer-based price prediction model for day-ahead bidding and a long short-term memory (LSTM)-dynamic programming hybrid real-time bidding.
    \item We build a transformer-based model for predicting electricity real-time prices and train it over the New York Independent System Operator (NYISO) market.
    \item We evaluated our proposed framework using actual price data obtained from four NYISO zones. The findings revealed that the model demonstrated a notable improvement, yielding up to a 29\% increase in profit when compared to solely participating in the RTM.
\end{itemize}
Notably, to our knowledge, this will be the first study in the literature to demonstrate the benefit of systematic day-ahead and real-time bidding using real-world price data. Results from this paper will empower storage participants to systematically bid into two-settlement markets, generating more profits and social welfare.

\subsection{Related Work}
Electricity price forecasting has been investigated heavily in the past years~\cite{Weron2,Lago}. Electricity is still hard to store and demand must be instantaneously balanced with demand in power systems, making the price of electricity have unique dynamics as it has special shapes for each time level and even sudden peaks. Many models have been developed previously to solve this problem; each has its own advantages and disadvantages~\cite{Skantze}. These models vary from statistical models~\cite{Deng}, economic equilibrium models~\cite{Hobbs}, agent-based models~\cite{Visudhiphan}, to experimental models~\cite{Denton}. Generally, there are three main types of price forecasting, short-term, mid-term, and long-term, focusing on forecasting hours ahead, weeks to months ahead, and years ahead. The horizon changes based on the needed application. This study focuses on the short-term horizon as it strongly influences the day-ahead and real-time markets operations.

Artificial intelligence (AI) models have risen as a leading solution, overcoming many limitations inherent in alternative approaches. Models such as multi-layered perceptron (MLP) have been shown to provide satisfactory performance and adaptability~\cite{PD16}. Furthermore, LSTM models have emerged as particularly suitable for price forecasting, outperforming hybrid methodologies like convolutional neural networks LSTM (CNN-LSTM) ~\cite{LSH22}. A comprehensive review of AI and statistical models in electricity price forecasting, encompassing datasets from diverse regions, is provided by~\cite{LMD21}.

Energy storage arbitrage in electricity wholesale markets has experienced rapid growth in recent years \cite{mcgrath2022battery}. Storage entities in wholesale electricity markets can participate in arbitrage by charging during periods of low prices and discharging during periods of high prices, thereby maximizing their profits. To evaluate potential profits, various models have been introduced in the literature, including price taker and strategic-behavior models \cite{sioshansi2021energy}. Price taker models operate under the assumption that their actions do not influence market operations, hence, they take the same price offered by the market. This simplifies the optimization problem by disregarding power balance constraints. Several approaches have been proposed to address this challenge, with Model Predictive Control (MPC) models often leading the literature due to their modeling flexibility \cite{arnold,abdulla2016optimal}. Additionally, SDP~\cite{Zheng} and RL~\cite{wang} models has been proposed frequently.

Two-settlement markets offer various participation options, including participation in DAM, RTM, and virtual bidding (VB). DAM participants are required to submit their bids one day in advance, necessitating accurate price forecasting for 24 hours-ahead. \cite{hashmi2020arbitrage} proposes an MPC model using forecasted prices. On the other hand, RTM participants submit bids for 5-minute intervals during operational hours, facing highly volatile prices. \cite{cao2020deep} introduces a deep reinforcement learning (DRL) model aimed at discovering an optimal policy to maximize arbitrage profits in RTM. \cite{Zheng} propose an analytical SDP approach to address real-time price arbitrage with efficient computation . Additionally, \cite{yousuf} outline a highly transferable price arbitrage model that combines model-based dynamic programming with neural networks.

VB option offers market participants the flexibility to place bids without the need for physical assets, enabling them to sell in the DAM and buy the same amount in the RTM, or vice versa. This option aims to improve market efficiency by minimizing the disparity between DAM and RTM prices \cite{li2015efficiency}. Additionally, \cite{li2021machine} discusses a virtual bidding framework utilizing a gradient boosting tree-based algorithm that takes price sensitivity into account. Market participants have the opportunity to participate in both DAM and RTM simultaneously, aiming to enhance their overall arbitrage profits. However, navigating multi-stage decision-making poses challenges, given the uncertainties in electricity prices. The proposed framework by \cite{rad2015} outlines an optimal bidding strategy for energy storage arbitrage across DAM and RTM, albeit without factoring in price uncertainty. Furthermore, \cite{krishnamurthy2017energy} have introduced an SDP model for storage arbitrage in DAM and RTM using conventional statistical methods to model the uncertainties.

The remainder of the paper is organized as follows: Section II introduces the methodology. Section III describes the dataset, and Section IV presents the case study. Section V concludes the paper.

\section{Energy Storage Arbitrage Framework}
\begin{figure}[htbp]
\centerline{\includegraphics[width=.9\linewidth]{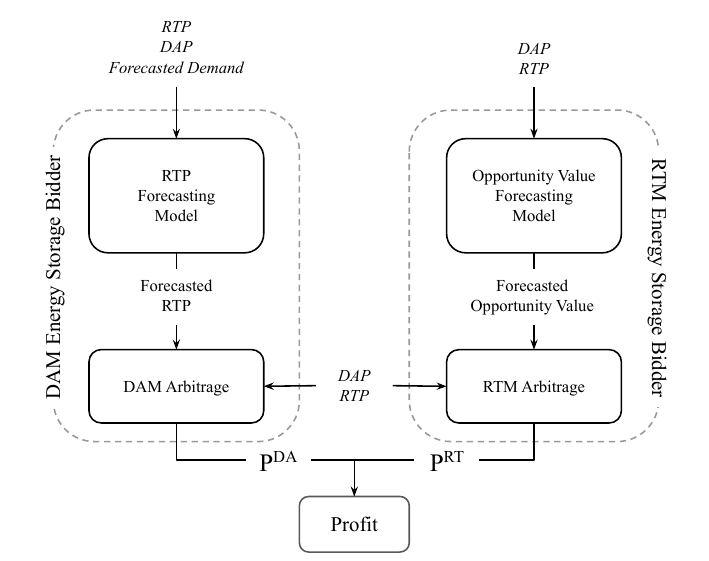}}
\caption{Pipeline of the proposed solution.}
\label{Full_system_diagram}
\end{figure}
Our profit-maximizing framework consists of a DAM bidding stage and sequential RTM bidding stages on a daily basis.
The {DAM stage} uses the predicted hourly average real-time price as the charge and discharge bid for the same hour in DAM. The {RTM stage} designs real-time bids sequentially for each hour to maximize real-time arbitrage profits.
Fig. \ref{Full_system_diagram} shows our full framework. In this framework, the DAM bids are based on predicted RTM prices, but the RTM bids are independent of DAM settlements. In the remainder of this section, we will present the theoretical analysis for the presented framework, our price prediction model, and the real-time arbitrage model.

\subsection{Problem Statement}
We assume the energy storage is a price-taker and submits economic bids to both DAM and RTM. The storage does not model the impact of its action over price formations and hence does not seek to exercise market power. Each bid, encompassing both a charging (purchase) and discharging (sale) bid, is submitted for every market-clearing hour. Each day comprises 48 bids: 24 for the DAM submitted collectively before the operating day starts and 24 for real-time, where bids are placed hourly with clearances at five-minute intervals. We thus formulate the two-settlement arbitrage problem as follows:
\begin{subequations}\label{ps}
\begin{align}
\max_{B\in\mathcal{B}} &\; \mathbb{E}_{\bm{\lambda}, \bm{\pi}\in \mathcal{M}}\Big\{\sum_{t\in\mathcal{T}}\lambda_t(p_t\up{D}-b_t\up{D}) - \sum_{s\in \mathcal{S}} c p_s\up{R}\nonumber\\
& +\sum_{s\in \mathcal{S}}\pi_s (p_s\up{R}-b_s\up{R} - {p_{\tau|s\in\tau}\up{D}}/{S}-{b_{\tau|s\in\tau}\up{D}}/{S})\Big\} \label{ps:obj}\\
&\{p_t\up{D},b_t\up{D},p_s\up{R},b_s\up{R}\} \in B \;\forall t\in\mathcal{T}, s\in\mathcal{S} \label{ps:c7}
\end{align}
in which
\begin{itemize}
    \item $t$ is the index of the hour at which the day-ahead market clears, the set $\mathcal{T}$ represents all hours from one day.
    \item $s$ is the index of real-time market clearing intervals, usually per five minutes, the set $\mathcal{S}$ represents all real-time market intervals in one day.
    \item $S$ is the number of real-time market clearing intervals in one day-ahead market interval. In our case study, we consider hourly day-ahead market clearing and 5-minute real-time market clearing intervals, so $S$ will be 12.
    \item $\lambda_t$ is the day-ahead market clearing price.
    \item $\pi_s$ is the real-time market clearing price.
    \item $p_t\up{D}$ and $b_t\up{D}$ are the storage's discharging and charging energy in day-ahead clearing results, respectively;
    \item $p_s\up{R}$ and $b_s\up{R}$ are the storage's \emph{physical} discharging and charging energy in real-time, respectively;
    \item $c$ is the cost to discharge; note that we assume the physical cost, such as degradation, only occurs during the discharge stage. This is to align with conventional generation cost models.
    \item $B$ is the bidding policy while $\mathcal{B}$ represents the set of all feasible bidding policies.
\end{itemize}
The objective function maximizes the arbitrage profit in DAM and RTM under the expectation that the day-ahead price and real-time price are a stochastic process to the market clearing model ($\bm{\lambda}, \bm{\pi}\in \mathcal{M}$). The first term represents DAM revenue, the second term is the physical discharge cost, and the third term represents real-time re-settlement based on day-ahead commitments. Note the subscript $\tau|s\in\tau$ represents the hour corresponding to the current RTM interval. 

In this objective, if we set $p_t\up{D}$ and $b_t\up{D}$ to zero, then the first term becomes zero and formulation becomes a real-time arbitrage objective. On the other hand, if we set $p_s\up{D} = p_{\tau|s\in\tau}\up{D}/{S}$ and $b_t\up{D}= b_{\tau|s\in\tau}\up{D}/{S}$, then it indicates the storage's physical charge and discharge power fully follow the DAM clearing results. Hence, the third term becomes zero and the storage is only participating in the DAM. \eqref{ps:c7} requires all storage actions must be the result of a bidding policy $B$, while $B$ must come from the set of all feasible bidding policy $\mathcal{B}$, which are policies that satisfy market clearing procedures and the following physical constraints ($\forall t\in\mathcal{T}, s\in\mathcal{S}$):
\begin{gather}
0\leq p_t\up{D},b_t\up{D} \leq P \label{ps:c1}\\
0\leq p_t\up{R},b_t\up{R} \leq P/S \label{ps:c2}\\
e_t\up{D}-e_{t-1}\up{D} = -p_t\up{D}/\eta + b_t\up{D}\eta \label{ps:c3}\\
e_s\up{R}-e_{s-1}\up{R} = -p_s\up{R}/\eta + b_s\up{R}\eta \label{ps:c4}\\
0\leq e_t\up{D} \leq E \label{ps:c5}\\
0\leq e_s\up{R} \leq E \label{ps:c6}\\
p_s^{R}=0\;\mathrm{if}\; \pi_t < 0 \label{ps:c8} 
\end{gather}
\end{subequations}
where $P$ is the storage power rating, $\eta$ is the efficiency, $e_t\up{D}$ is the state-of-charge (SoC) constraint in DAM clearing, while $e_s\up{R}$ is the physical SoC constraint. \eqref{ps:c1} and \eqref{ps:c2} model the power rating constraint in day-ahead and real-time operations, \eqref{ps:c3}, \eqref{ps:c4}, \eqref{ps:c5} and \eqref{ps:c6} model the SoC evolution constraints in the same manner. \eqref{ps:c8} ensures the storage doesn't discharge when the price is negative, a sufficient condition to avoid simultaneous charging and discharging~\cite{yousuf}. Note that this constraint is only necessary in RTM as day-ahead prices are seldom negative.

\subsection{Decomposition Reformulation }

Using \eqref{ps} to directly find the optimal arbitrage policy can be extremely difficult as the electricity prices follow highly stochastic multi-stage processes and the bidding requirements are also non-trivial to model. Instead, we temporally relax the bidding constraint \eqref{ps:c7} and treat the powers as direct decision variables. The goal thus is to establish intuition in the market participation and use it to design bids. We reformulate \eqref{ps:obj} by partially merging the first and the third term to decouple day-ahead and real-time decision variables and stochastic processes, the resulting formulation is:
\begin{align}
    & \textstyle\max_{p_t\up{D}, b_t\up{D}} \; \mathbb{E}_{\bm{\lambda}} \Big\{ \sum_{t\in\mathcal{T}}\lambda_t(p_t\up{D}-b_t\up{D}) - \overline{\pi}_t(p_t\up{D}-b_t\up{D}) \Big\}\nonumber\\
& \textstyle + \max_{p_s\up{R}, b_s\up{R}}\mathbb{E}_{\bm{\pi}}\Big\{\sum_{s\in \mathcal{S}}\pi_s (p_s\up{R}-b_s\up{R} ) - \sum_{s\in \mathcal{S}} c p_s\up{R}\Big\}\label{ref}
\end{align}
where $\overline{\pi}_t$ is the expected hourly average real-time price.
This reformulation is based on the following approximation to decouple price expectation from power expectation.
\begin{proposition}\emph{Real-time price expectation under price-taker assumptions.}
The DAM clearing results are statistically independent of real-time price realizations. This allows us to convert the whole-day summation of RTM intervals ($s\in \mathcal{S}$) into summation of hours ($t\in \mathcal{T}$) and summation of intervals within a given hour ($s\in t$), hence
    \begin{align}
        \mathbb{E}_{\bm{\pi}}\Big\{\sum_{s\in\mathcal{S}}{\pi_s p_{\tau|s\in\tau}\up{D}}/{S}\Big\} & = \sum_{t\in \mathcal{T}} \mathbb{E}_{\bm{\pi}}\Big\{\sum_{s\in t} \pi_s {p_{\tau|s\in\tau}\up{D}}/{S}\Big\} \nonumber \\
       & = \sum_{t\in \mathcal{T}}(\mathbb{E}_{\bm{\pi}}\sum_{s\in t}\pi_s)p_t\up{D} = \sum_{t\in \mathcal{T}}\overline{\pi}_tp_t\up{D}\nonumber
    \end{align}
\end{proposition}
\noindent and similarly to $b_t\up{D}$.

\begin{proof}
This proposition is based on the price-taker assumption that the storage does not assume its action would impact the market price and will not try to exercise market power, then it is trivial that $p_{\tau}\up{D}$ will not change the distribution of price $\pi_s$. This is also a common practice in bidding designs~\cite{Bitar, xu2018optimal} when assuming participants are price takers in the market
\end{proof}

Based on the previous proposition, we state our second proposition
\begin{proposition}\emph{Two-stage reformulation.}
    \eqref{ps:obj} and \eqref{ref} are equivalent formulations under the price-taker assumption.
\end{proposition}

\begin{proof}
    Proposition~1 allows us to decouple the price expectation from the action, and further allows us to decouple the stochastic model of $\bm{\lambda}$ and $\bm{\pi}$. Thus, we only need to justify the separation of the maximization objective from \eqref{ps:obj} to \eqref{ref}, which is intuitively given there is no coupling physical constraint between DAM and RTM decisions. Thus, \eqref{ps:obj} and \eqref{ref} are equivalent.
\end{proof} 

By observing the format of the reformulated objective \eqref{ref}: the first maximization term (associated with constraints \eqref{ps:c1}, \eqref{ps:c3}, and \eqref{ps:c5}) corresponds to the day-ahead bidding problem with the bidding cost equal to the expected average real-time price over the corresponding hour; the second maximization term (associated with constraints \eqref{ps:c2}, \eqref{ps:c4}, \eqref{ps:c6} and \eqref{ps:c8}) is a real-time arbitrage model without considering any day-ahead results. Thus, this reformulation justifies our bidding model as shown in Figure~\ref{Full_system_diagram}, which also satisfies market bidding requirements as in the original formulation \eqref{ps}.

\subsection{Electricity Price Forecasting}

According to our decomposition reformulation result, the key to designing day-ahead bids is to predict the real-time price of the same hour. Thus in this subsection, we will present a machine learning model for predicting real-time prices.

\begin{figure}[htbp]
\centerline{\includegraphics[width=\linewidth]{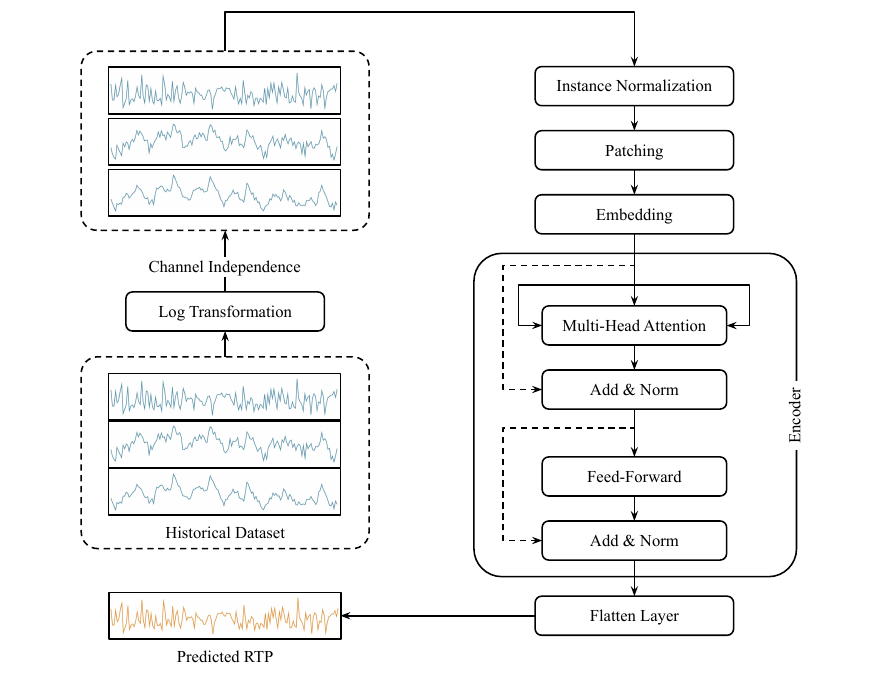}}
\caption{Our real-time price forecasting model overview.}
\label{our_model}
\end{figure}

In the first stage of our framework, we forecast the real-time price, which serves as a critical input for our day-ahead arbitrage model. The foundation of our forecasting model is built upon the patch time series transformer (PatchTST), a model primarily developed as outlined by \cite{Nie, zhou2021informer}.
Transformers~\cite{Vaswani}, a recent class of deep learning models initially designed for applications in natural language processing, have consistently demonstrated remarkable performance across various domains. Among these, the PatchTST~\cite{Nie} stands out as a prominent model for time-series forecasting, consistently surpassing other transformers. Fig.~\ref{our_model} shows an overview of our real-time price forecasting model. It only uses the encoder from the traditional transform and has two main features: channel-independence, split the multivariate time series into univariate time series that share the same weights, and each channel may have a different attention map; patching, decompose the time series into subsequences, capture local information, and reduce model complexity.

We follow a similar formulation to \cite{Nie}, but with the input to be a multivariate time-series with a look-back window $L$ and the output a single univariate time-series with a forecasting horizon of $T$ as described by the following equations:

\begin{equation}
 \boldsymbol{x} =
\begin{bmatrix}
x_{1,1} & \cdots & x_{1,L} \\
\vdots  & \ddots & \vdots  \\
x_{M,1} & \cdots & x_{M,L} \\
\end{bmatrix}, ~~~
\boldsymbol{\hat{x}} = \{\hat{x}_{1,1}, \ldots, \hat{x}_{1,T}\}
\label{eq:output}
\end{equation}

The process begins by breaking down the multivariate time-series input into several distinct univariate time-series during a phase called channel-independence. Each of these individual univariate time series is then separately fed into the transformer backbone.

The transformer backbone primarily comprises three key components: instance normalization, patching, and the traditional transformer encoder. Every univariate time series will first be normalized with zero mean and unit standard deviation and then divided into a sequence of patches. This patching not only minimizes the number of input tokens required for the encoder but also serves to reduce computational complexity. Moreover, it helps the model to effectively capture extended historical sequence, enhancing its predictions.

The primary purpose of the transformer encoder is to transform input patches into a different representation. This transformation involves several steps. Initially, a linear projection of the input is obtained. Subsequently, scaled dot-product attention is applied, as expressed by the equation:
\begin{align}
\text{Attention}(Q, K, V) &= \textstyle\text{softmax}\left(\frac{Q \cdot K^T} {\sqrt{d_k}}\right) V 
\label{eq:scaled}
\end{align}

Here, the matrices $Q$, $K$, and $V$ correspond to the query, key, and value matrices, respectively, and have been transformed by each head in the multi-head attention mechanism. Following that, the resulting attention values are batch-normalized and passed through a feed-forward neural network before reaching the flattening layer, ultimately yielding the predictions as shown in~\eqref{eq:output}.

Finally, the loss function that our model minimizes during training is the mean square error (MSE), which we employ to guide the model's optimization process. This loss function is mathematically defined as follows:
\begin{align}
\text{MSE} = \textstyle\frac{1}{T} \sum_{t=1}^{T} (x_t - \hat{x}_t)^2
\label{eq:loss}
\end{align}

\subsection{Real-Time Energy Storage Arbitrage Model}

Our storage arbitrage model for real-time bidding is based on~\cite{yousuf}, which solves the real-time arbitrage problem (second part of \eqref{ref}) following a non-anticipatory bidding policy. The model predicts the opportunity value of the state of charge (SOC) and then maximizes the storage arbitrage profit. The model combines model-based dynamic programming with neural networks. Unlike MPC, which is very sensitive to price volatility over a 24-hour look-ahead horizon, this method focuses on predicting the value function within a single time period. Furthermore,~\cite{yousuf} is highly transferable, which makes it suitable for different markets and participation scenarios. It has been evaluated using NYISO price data; it achieved a profit ratio ranging from 70\% $\sim$ 90\% when compared to the perfect forecast. These features make this method ideal for our system. Due to page limitations, we will not present the full formulation in this paper as this model has already been published.

\section{Data Description}

This work utilizes a dataset spanning a five-year period, ranging from 2017 to 2021. This dataset encompasses multiple price zones within the NYISO market, including NYC, LONGIL, NORTH, and WEST:

\begin{itemize}
    \item Hourly Day-Ahead Price (DAP): This variable represents the smallest available resolution for day-ahead LMP. It has been used in all parts of the framework, and it has a direct influence on both day-ahead and real-time markets.
    \item 5-minute Real-Time Price (RTP): This is the smallest resolution in NYISO. These granular data points help in finding the best solution by the real-time arbitrage model.
    \item Hourly Real-Time Price (RTP): This resolution was added to match the resolution of the day-ahead prices when used for the day-ahead arbitrage model.
    \item Hourly demand forecasting: This variable provides insights into the predicted energy demand on an hourly basis. It's critical to understand how real-time pricing is influenced by changes in energy demand.
\end{itemize}

\begin{figure}[htbp]
\centerline{\includegraphics[width=\linewidth]{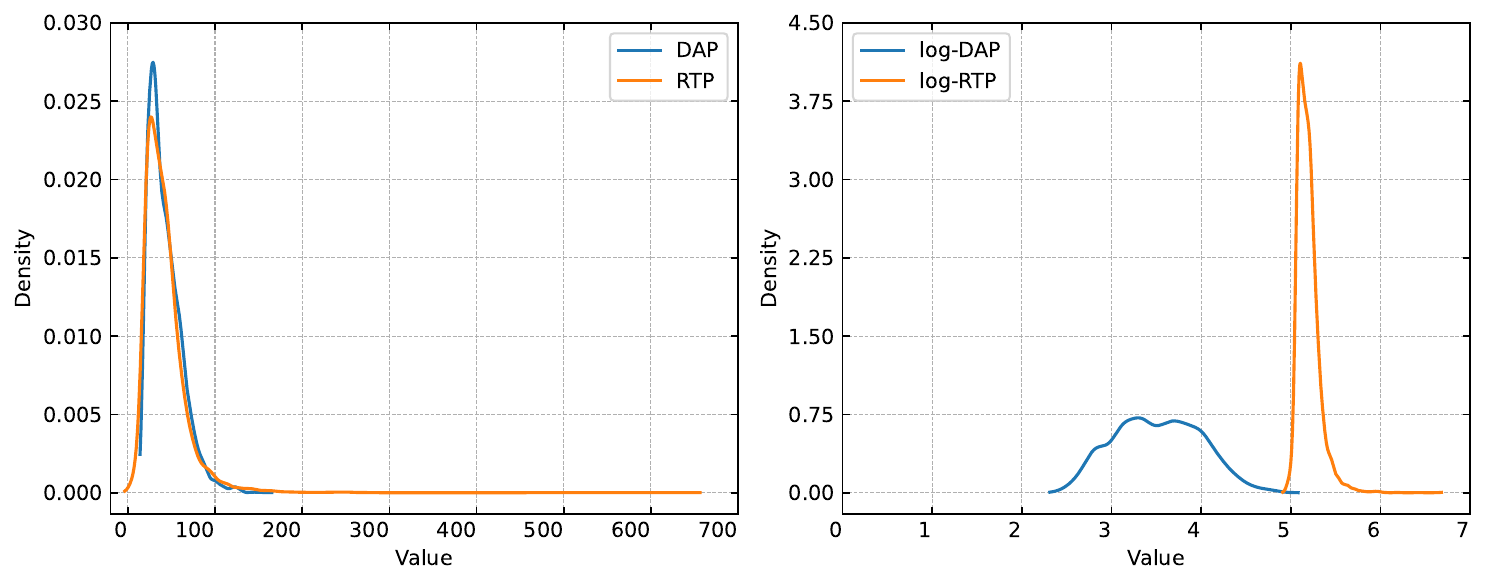}}
\caption{Density distribution of the day-ahead and real-time prices for NYC Zone price during 2021.}
\label{fig:NYC_density}
\end{figure}

Our framework is constructed around these variables, driven by their direct impact on DAM and RTM, as well as their availability. These variables encapsulate pricing dynamics from both DAM and RTM. Fig.~\ref{fig:NYC_density} shows the density distribution of the prices in the testing set. We can see that real-time prices have more outliers than the day-ahead prices. Also, the log-transformed prices gave a better distribution, hence expected to improve the forecasting task. Table \ref{table:zones_stats} shows statistics about the historical price data for the four NYISO zones, it is mainly the price mean ($\mu$), the price standard deviation ($\sigma$), and the number of negative real-time prices.
\begin{table}[htbp]
    \centering
    \caption{NYISO Price Zones's Statistics}
    \begin{tabularx}{\linewidth}{>{}X> {\centering\arraybackslash}X>{\centering\arraybackslash}X>{\centering\arraybackslash}X>{\centering\arraybackslash}X>{\centering\arraybackslash}X}
        \toprule
        \textbf{Zone} & \textbf{$\mu$ DAP} & \textbf{$\sigma$ DAP} & \textbf{$\mu$ RTP} & \textbf{$\sigma$ RTP} & \textbf{Negative RTP} \\
        \midrule
        \textbf{NYC} & 42.54 & 18.94 & 42.49 & 26.89 & 1 \\
        \textbf{LONGIL} & 54.11 & 31.40 & 54.85 & 63.72 & 42 \\
        \textbf{NORTH} & 22.65 & 16.07 & 23.51 & 30.27 & 937 \\
        \textbf{WEST} & 31.03 & 19.47 & 30.89 & 27.08 & 44 \\
        \bottomrule
    \end{tabularx}
    \label{table:zones_stats}
\end{table}

\section{Experiment Setup and Results}
The study was performed using data collected from NYISO over the span of five years (2017-2021), four years for training and validation and one year for testing. Real-time prices, day-ahead prices, and demand forecasting were chosen as features due to their direct influence on the real-time price. A look-back window of two weeks and the prediction horizon of 24 hour-ahead were considered. The data were normalized before being fed into the model and the price features were log-transformed to reduce the impact of the extreme spikes on the model predictions. We consider battery storage with a 2-hour storage duration, 1 MWh storage capacity, 90\% one-way efficiency, and \$10/MWh fixed discharge cost that covers degradation cost and other operational costs.

All experiments are run on a high performance computing cluster with Intel Xeon Platinum 8640Y 2Ghz CPU, 512GB Memory, and 2xNVIDIA A40 GPU. We used Tenserflow 2.15.0, CVXPY 1.4.2, PyTorch 2.2.0, and CUDA 12.1. All codes are available on Github\footnote{https://github.com/Alghumayjan/Energy-Storage-Arbitrage-in-Two-Settlement-Markets}.

\subsection{Electricity Prediction Results}
In this subsection, we present the results of the electricity real-time price forecasting model. We benchmark the model with the DLinear model~\cite{zeng}, which was shown to beat most time-series transformers.
The DLinear is a special linear model designed for time-series problems; it utilizes the decomposition layer of the Autuoformer model~\cite{wu}, where it splits the input data into multiple components that represent the trend and seasonality in the input and then predicts the output based on these components. Table \ref{table:models_mse} shows a comparison of our model (the transformer) and the baseline model (the DLinrear) in terms of the prediction accuracy on the testing set using MSE for the NYC price zone. We can see that the our model beats the baseline model. Also, we can see that the model with log transformation achieved lower errors in general.
\begin{table}[ht]
    \caption{Testing Set Prediction Accuracy}
    \centering
    \renewcommand{\arraystretch}{1.5} 
    \begin{tabular}{c c c c c}
        \toprule
        \textbf{Model} & \textbf{Ours} & \textbf{log-Ours} & \textbf{Baseline} & \textbf{log-Baseline} \\
        \midrule
        \textbf{MSE} & 0.458 & \textbf{0.429} & 0.506 & 0.457 \\
        \bottomrule
    \end{tabular}
    \label{table:models_mse}
\end{table}

We have considered different metrics to evaluate the models. Based on~\cite{Lago}, one of the most suitable metrics for our case is the relative mean absolute error (rMAE)~\cite{hyndman}, which is the ratio of the MAE resulting from the original model and some naive model. \eqref{eq:rMAE} shows how the rMAE is calculated where the variable $x$ represents the actual output values, while $\hat{x}$ represents the forecasted output values. The naive model is described as $\hat{x}_{d,h}^{\mathrm{naive}} = x_{d-7,h}$ where in our case, we consider the 7-day lag model, where the forecast is the same hour one week ago. Table~\ref{table:models_mae} shows the results based on the MAE and rMAE metrics after rescaling the predictions to the actual levels. These results align with the results from MSE, where our model with log tranformation has the lowest error.
\begin{align}
\text{rMAE} = \frac{\mathrm{MAE}}{\mathrm{MAE}_{\mathrm{naive}}} =  \frac{\sum_{d=1}^{365}\sum_{h=1}^{24} |x_{d,h} - \hat{x}_{d,h}|}{\sum_{d=1}^{365}\sum_{h=1}^{24} |x_{d,h} - \hat{x}_{d,h}^{\mathrm{naive}}|}
\label{eq:rMAE}
\end{align}
\begin{table}[ht]
    \caption{Rescaled Testing Set Prediction Accuracy}
    \centering
    \renewcommand{\arraystretch}{1.5} 
    \begin{tabular}{c c c c c}
        \toprule
        \textbf{Model} & \textbf{Ours} & \textbf{log-Ours} & \textbf{Baseline} & \textbf{log-Baseline} \\
        \midrule
        \textbf{MAE} & 11.01 & \textbf{10.62} & 11.95 & 10.97 \\
        \textbf{rMAE} & 0.62 & \textbf{0.60} & 0.67 & 0.62 \\
        \bottomrule
    \end{tabular}
    \label{table:models_mae}
\end{table}
\begin{table*}[h]
\centering
\caption{Model Performance Summary}
\label{tab:model_performance}
\begin{tabular*}{\textwidth}{@{\extracolsep{\fill}}lccccccccc}
\toprule
\multirow{2}{*}{Price Zone} & \multirow{2}{*}{RTP mean} & \multicolumn{4}{c}{log-Ours} & \multicolumn{4}{c}{log-Baseline} \\
\cmidrule(lr){3-6} \cmidrule(lr){7-10}
& & MAE (\$/MWh) & rMAE & Profit (\$) & IPM & MAE (\$/MWh) & rMAE & Profit (\$) & IPM \\
\midrule
NYC & 42.74 & 10.62 & 0.60 & \textbf{13875} & 29\% & 10.97 & 0.62 & 12580 & 18\% \\
LONGIL & 55.05 & 27.08 & 0.76 & \underline{33849} & 8\% & 29.40 & 0.82 & 32786 & 5\% \\
NORTH & 23.53 & 13.67 & 0.68 & \textbf{20246} & 18\% & 17.87 & 0.89 & 17318 & 1\% \\
WEST & 31.03 & 12.00 & 0.68 & \textbf{16836} & 26\% & 13.59 & 0.77 & 15289 & 14\% \\
\bottomrule
\end{tabular*}
\end{table*}
\begin{table*}[h]
\centering
\begin{tabular*}{\textwidth}{@{\extracolsep{\fill}}lccccccccc}
\multirow{2}{*}{Price Zone} & \multirow{2}{*}{RTP mean} & \multicolumn{4}{c}{Ours} & \multicolumn{4}{c}{Baseline} \\
\cmidrule(lr){3-6} \cmidrule(lr){7-10}
& & MAE (\$/MWh) & rMAE & Profit (\$) & IPM & MAE (\$/MWh) & rMAE & Profit (\$) & IPM \\
\midrule
NYC & 42.74 & 11.01 & 0.62 & \underline{13802} & 29\% & 11.95 & 0.67 & 12564 & 18\% \\
LONGIL & 55.05 & 26.06 & 0.73 & \textbf{33995} & 9\% & 25.95 & 0.73 & 33727 & 8\% \\
NORTH & 23.53 & 13.41 & 0.66 & \underline{19490} & 14\% & 14.84 & 0.74 & 17856 & 4\% \\
WEST & 31.03 & 11.70 & 0.66 & \underline{16785} & 26\% & 12.19 & 0.69 & 15345 & 15\% \\
\bottomrule
\end{tabular*}
\end{table*}

\begin{figure*}[htbp]
\centerline{\includegraphics[width=\linewidth]{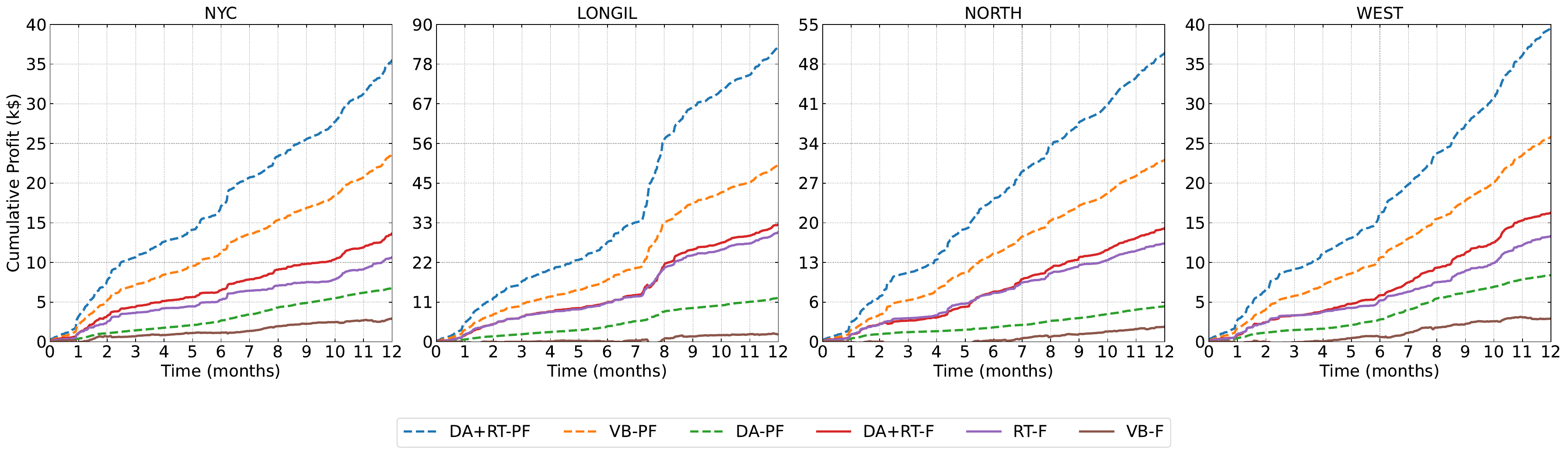}}
\caption{The results of the system for all NYISO four price zones during 2021. Note that the dashed lines are participating options using perfect real-time price forecast, and the solid lines are participating options using forecasted real-time price generated by the log-Transformer model.}
\label{fig:arb_profit}
\end{figure*}
\subsection{Combined Profit Results}
We now present the results of the arbitrage simulation for the entire proposed framework described in Section II. In this study, we investigate multiple settings and participating options. It mainly participates in DAM only, RTM only, VB, or both DAM and RTM. We consider these options under the perfect forecast real-time price case and the forecasted real-time price case:
\begin{itemize}
    \item \textbf{DA+RT-PF, DA-PF, VB-PF}: participates in both DAM and RTM, only in the DAM, or only in VB, respectively, using perfect forecast real-time price.
    \item \textbf{DA+RT-F, RT-F, VB-F}: participates in both DAM and RTM, only in the DAM, or only in VB, respectively, using forecasted RT price via the developed models.
\end{itemize}

Table \ref{tab:model_performance} presents a summary model's performance for all the models in terms of prediction accuracy and total profit for the testing year (2021) across the four NYISO price zones considered. The profit presented is for the combined market participation option (DA+RT-F). We also report the improved profit margin (IPM), which indicates the percentage of the increased profit we got by participating in the DAM and RTM over participating in the RTM only. 

The log-ours model has the highest (Bold) or second-highest (underlined) profit across all zones. In addition,  predictor errors do not strictly correlate with profits. For example, the model with the lowest error in the LONGIL zone is the baseline, but the model with the highest profit is ours. Moreover, the LONGIL zone has the highest prediction errors and the lowest IPM in all the models, which can be explained by Table \ref{table:zones_stats}; LONGIL has a very high standard deviation in the real-time price, making it highly volatile and harder to predict. Furthermore, note that the baseline models suffer a lot in the NORTH zone, as they have IPM lower than 4\% in comparison to our models, which have IPM higher than 14\%; this suggests that the baseline models is too sensitive to negative prices, while our models are more robust. In general, the results show that our models is more profitable and reliable.  

Fig. \ref{fig:arb_profit} shows the cumulative profit for all participating options considered in the study. We chose the log-ours model in all cases for its overall best performance. The proposed joint bidding approach provides the highest profit among all cases when assuming realistic price predictions. On the other hand, VB with perfect forecasts and DA+RT-PF case, which also includes VB, significantly surpasses all other approaches as it can fully take advantage of real-time price spikes. However, it is unrealistic to expect such high yields as price spikes are extremely difficult to predict~\cite{zhao2007framework, baltaoglu2018algorithmic}, and the VB profit becomes the lowest when using realistic price forecasts. 

\begin{figure}[htbp]
\centerline{\includegraphics[width=\linewidth]{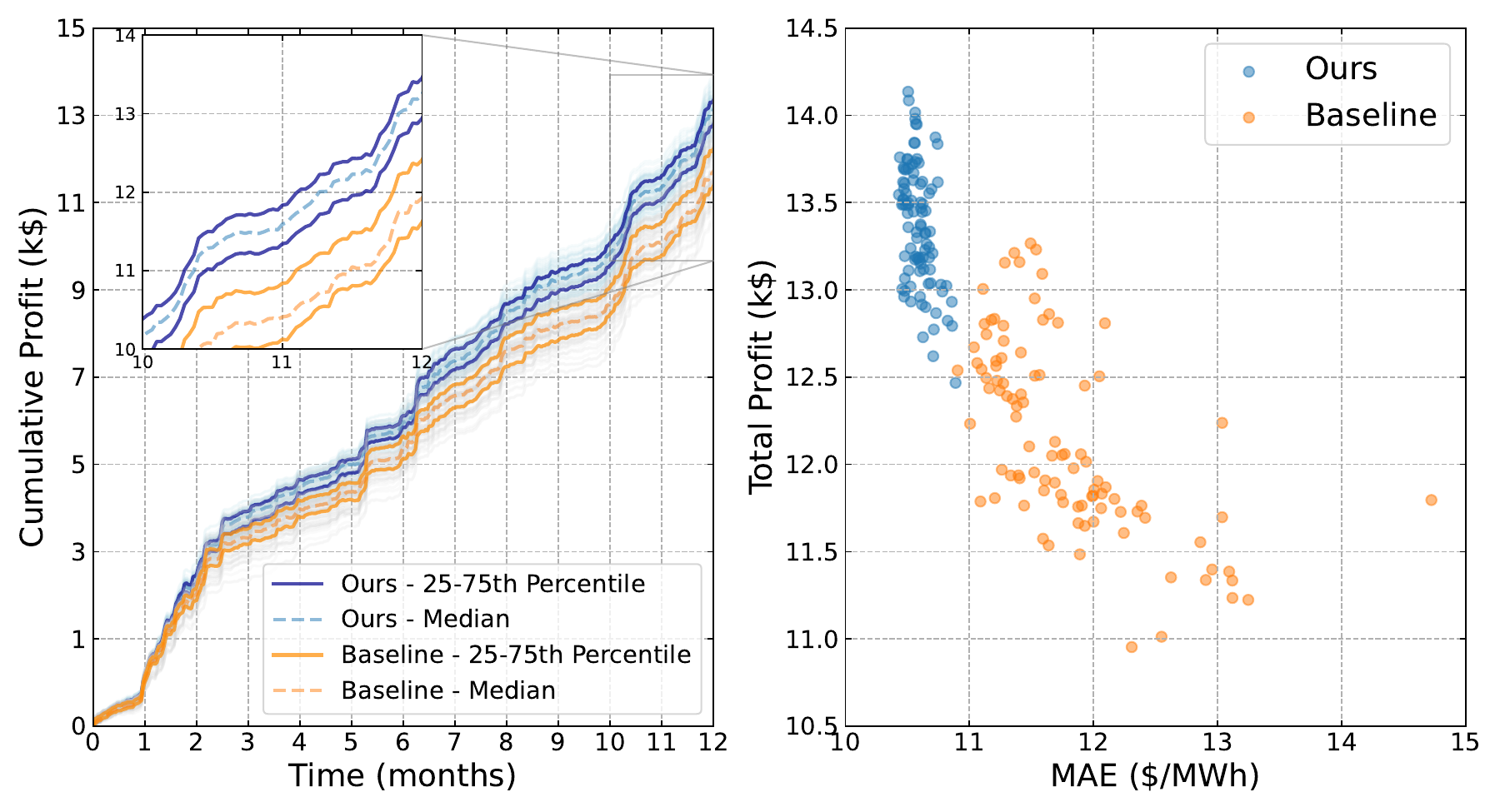}}
\caption{Left: The cumulative profit of the trained models with randomly sampled hyperparameters. Right: the error of the models and their profits.}
\label{fig:model_sensitivity}
\end{figure}
\subsection{Sensitivity Analysis}
To assess the robustness of our bidding strategy against errors in real-time price forecasting, we demonstrate how the profitability of our proposed approach degrades with the performance of the real-time price forecaster. We have trained multiple versions of the real-time price forecaster, each with hyperparameters sampled from a uniform distribution. Subsequently, we conducted arbitrage simulations using the resulting forecasted prices. The left plot in Fig. \ref{fig:model_sensitivity} illustrates the cumulative profit of our approach compared to the baseline. Meanwhile, the right plot shows the errors of the predictions against the total profit. The figure shows that lower prediction errors correspond to higher profits. Furthermore, our model consistently exhibits superior profitability and reduced errors in general.

\begin{table}[htbp]
    \centering
    \caption{Number of negative profit days of various participation}
    \begin{tabularx}{\linewidth}{>{\centering\arraybackslash}X>{\centering\arraybackslash}X>{\centering\arraybackslash}X>{\centering\arraybackslash}X>{\centering\arraybackslash}X}
        \toprule
        \textbf{Zone} & \textbf{VB} & \textbf{RT} & \textbf{RT+DA} & \textbf{IPM} \\
        \midrule
        NYC & 104 & 90 & 59 & 29\% \\
        LONGIL & 153 & 38 & 39 & 8\% \\
        NORTH & 111 & 28 & 27 & 18\% \\
        WEST & 115 & 39 & 33 & 26\% \\
        \bottomrule
    \end{tabularx}
    \label{t-v}
\end{table}
\subsection{Risk Analysis}
By analyzing the cumulative profit curves and daily profit distributions, we observe that when predicted prices are utilized, many instances result in negative profits, posing a potential threat to energy storage owners. Ideally, we aim for results that closely resemble the scenarios with perfect forecasts. Table~\ref{t-v} provides insights into the frequency of daily negative profit occurrences for the three market participation options. Our analysis reveals that participating in both markets, while yielding on average 20\% higher profits than RTM-only, does not increase the risk of profit loss. In particular, in NYC, joint participation reduced negative profit days by 34\% compared to RTM-only. This indicates that combined market participation mitigates risk and enhances overall profitability or, at worst, maintains higher profitability levels while managing risk. VB, on the other hand, provides the highest number of negative profit days. This result reiterates that while VB is promising with perfect foresight, it is very risky in practice. 

\section{Conclusion}
We introduced an integrated model for optimizing energy storage bidding in two-settlement electricity markets. Combining a transformer-based model for day-ahead bidding and an LSTM-dynamic programming hybrid model for real-time bidding, we have demonstrated the potential to significantly enhance profit margins in two-settlement electricity markets. Our study addressed the problem of effectively managing energy storage bids in volatile real-time prices. Our experimentation with actual price data from the NYISO market shows that our proposed framework outperforms the sole real-time bidding approach. Specifically, we observed a substantial increase in profit, with improvements of up to 29\% when compared to strategies focused solely on the real-time market. Our risk analysis over negative profit days also shows that the proposed joint bidding algorithm provides less or similar risk than real-time market participation. In summary, using our proposed algorithm, bidding storage in day-ahead and real-time markets provides higher profits than real-time markets without added risk.

For the next steps, we plan to improve our forecasting model with an enhanced customized loss function to reflect the loss that arose from the decision-making process. In addition, we want to implement risk control measurements to reduce the risk that comes from price uncertainties.

% \end{thebibliography}
\bibliographystyle{IEEEtran}
\bibliography{references}

% Generated by IEEEtran.bst, version: 1.14 (2015/08/26)
\begin{thebibliography}{10}
\providecommand{\url}[1]{#1}
\csname url@samestyle\endcsname
\providecommand{\newblock}{\relax}
\providecommand{\bibinfo}[2]{#2}
\providecommand{\BIBentrySTDinterwordspacing}{\spaceskip=0pt\relax}
\providecommand{\BIBentryALTinterwordstretchfactor}{4}
\providecommand{\BIBentryALTinterwordspacing}{\spaceskip=\fontdimen2\font plus
\BIBentryALTinterwordstretchfactor\fontdimen3\font minus \fontdimen4\font\relax}
\providecommand{\BIBforeignlanguage}[2]{{%
\expandafter\ifx\csname l@#1\endcsname\relax
\typeout{** WARNING: IEEEtran.bst: No hyphenation pattern has been}%
\typeout{** loaded for the language `#1'. Using the pattern for}%
\typeout{** the default language instead.}%
\else
\language=\csname l@#1\endcsname
\fi
#2}}
\providecommand{\BIBdecl}{\relax}
\BIBdecl

\bibitem{ferc}
\text{Federal Energy Regulatory Commission}, ``Electric storage participation in markets operated by regional transmission organizations and independent system operators,'' \emph{Order No. 841, 162 FERC}, vol.~61, p. 127, 2018.

\bibitem{sioshansi}
R.~Sioshansi, P.~Denholm, and T.~Jenkin, ``Market and policy barriers to deployment of energy storage,'' \emph{Economics of Energy \& Environmental Policy}, vol.~1, no.~2, pp. 47--64, 2012.

\bibitem{Williams}
O.~Williams and R.~Green, ``Electricity storage and market power,'' \emph{Energy policy}, vol. 164, p. 112872, 2022.

\bibitem{us2021form}
\text{US Energy Information Association}, ``Form eia-860 detailed data with previous form data (eia-860a/860b),'' 2022.

\bibitem{arnold}
M.~Arnold and G.~Andersson, ``Model predictive control of energy storage including uncertain forecasts,'' in \emph{Power systems computation conference (PSCC), Stockholm, Sweden}, vol.~23.\hskip 1em plus 0.5em minus 0.4em\relax Citeseer, 2011, pp. 24--29.

\bibitem{abdulla2016optimal}
K.~Abdulla, J.~De~Hoog, V.~Muenzel, F.~Suits, K.~Steer, A.~Wirth, and S.~Halgamuge, ``Optimal operation of energy storage systems considering forecasts and battery degradation,'' \emph{IEEE Transactions on Smart Grid}, vol.~9, no.~3, pp. 2086--2096, 2016.

\bibitem{Zheng}
N.~Zheng, J.~Jaworski, and B.~Xu, ``Arbitraging variable efficiency energy storage using analytical stochastic dynamic programming,'' \emph{IEEE Transactions on Power Systems}, vol.~37, no.~6, pp. 4785--4795, 2022.

\bibitem{wang}
H.~Wang and B.~Zhang, ``Energy storage arbitrage in real-time markets via reinforcement learning,'' in \emph{2018 IEEE Power \& Energy Society General Meeting (PESGM)}.\hskip 1em plus 0.5em minus 0.4em\relax IEEE, 2018, pp. 1--5.

\bibitem{Milstein}
I.~Milstein and A.~Tishler, ``Can price volatility enhance market power? the case of renewable technologies in competitive electricity markets,'' \emph{Resource and Energy Economics}, vol.~41, pp. 70--90, 2015.

\bibitem{qin2023role}
X.~Qin, B.~Xu, I.~Lestas, Y.~Guo, and H.~Sun, ``The role of electricity market design for energy storage in cost-efficient decarbonization,'' \emph{Joule}, 2023.

\bibitem{Byrne}
R.~H. Byrne, T.~A. Nguyen, D.~A. Copp, R.~J. Concepcion, B.~R. Chalamala, and I.~Gyuk, ``Opportunities for energy storage in caiso: Day-ahead and real-time market arbitrage,'' in \emph{2018 International Symposium on Power Electronics, Electrical Drives, Automation and Motion (SPEEDAM)}.\hskip 1em plus 0.5em minus 0.4em\relax IEEE, 2018, pp. 63--68.

\bibitem{Caiso}
\text{California Independent System Operator (CAISO)}, ``Special report on battery storage,'' 2023.

\bibitem{rad2015}
H.~Mohsenian-Rad, ``Optimal bidding, scheduling, and deployment of battery systems in california day-ahead energy market,'' \emph{IEEE Transactions on Power Systems}, vol.~31, no.~1, pp. 442--453, 2015.

\bibitem{Weron2}
R.~Weron, ``Electricity price forecasting: A review of the state-of-the-art with a look into the future,'' \emph{International journal of forecasting}, vol.~30, no.~4, pp. 1030--1081, 2014.

\bibitem{Lago}
J.~Lago, G.~Marcjasz, B.~De~Schutter, and R.~Weron, ``Forecasting day-ahead electricity prices: A review of state-of-the-art algorithms, best practices and an open-access benchmark,'' \emph{Applied Energy}, vol. 293, p. 116983, 2021.

\bibitem{Skantze}
P.~L. Skantze and M.~Ilic, \emph{Valuation, hedging and speculation in competitive electricity markets: a fundamental approach}.\hskip 1em plus 0.5em minus 0.4em\relax Springer Science \& Business Media, 2001, vol. 643.

\bibitem{Deng}
S.~Deng, \emph{Stochastic models of energy commodity prices and their applications: Mean-reversion with jumps and spikes}.\hskip 1em plus 0.5em minus 0.4em\relax University of California Energy Institute Berkeley, 2000.

\bibitem{Hobbs}
B.~F. Hobbs, C.~B. Metzler, and J.-S. Pang, ``Strategic gaming analysis for electric power systems: An mpec approach,'' \emph{IEEE transactions on power systems}, vol.~15, no.~2, pp. 638--645, 2000.

\bibitem{Visudhiphan}
P.~Visudhiphan and M.~D. Ilic, ``Dynamic games-based modeling of electricity markets,'' in \emph{IEEE Power Engineering Society. 1999 Winter Meeting (Cat. No. 99CH36233)}, vol.~1.\hskip 1em plus 0.5em minus 0.4em\relax IEEE, 1999, pp. 274--281.

\bibitem{Denton}
M.~J. Denton, S.~J. Rassenti, V.~L. Smith, and S.~R. Backerman, ``Market power in a deregulated electrical industry,'' \emph{Decision Support Systems}, vol.~30, no.~3, pp. 357--381, 2001.

\bibitem{PD16}
I.~P. Panapakidis and A.~S. Dagoumas, ``Day-ahead electricity price forecasting via the application of artificial neural network based models,'' \emph{Applied Energy}, vol. 172, pp. 132--151, 2016.

\bibitem{LSH22}
M.~Lehna, F.~Scheller, and H.~Herwartz, ``Forecasting day-ahead electricity prices: A comparison of time series and neural network models taking external regressors into account,'' \emph{Energy Economics}, vol. 106, p. 105742, 2022.

\bibitem{LMD21}
J.~Lago, G.~Marcjasz, B.~De~Schutter, and R.~Weron, ``Forecasting day-ahead electricity prices: A review of state-of-the-art algorithms, best practices and an open-access benchmark,'' \emph{Applied Energy}, vol. 293, p. 116983, 2021.

\bibitem{mcgrath2022battery}
G.~McGrath and O.~Comstock, ``Battery systems on the us power grid are increasingly used to respond to price,'' \emph{Today in Energy}, 2022.

\bibitem{sioshansi2021energy}
R.~Sioshansi, P.~Denholm, J.~Arteaga, S.~Awara, S.~Bhattacharjee, A.~Botterud, W.~Cole, A.~Cortes, A.~De~Queiroz, J.~DeCarolis \emph{et~al.}, ``Energy-storage modeling: State-of-the-art and future research directions,'' \emph{IEEE transactions on power systems}, vol.~37, no.~2, pp. 860--875, 2021.

\bibitem{hashmi2020arbitrage}
M.~U. Hashmi, D.~Deka, A.~Bu{\v{s}}i{\'c}, L.~Pereira, and S.~Backhaus, ``Arbitrage with power factor correction using energy storage,'' \emph{IEEE Transactions on Power Systems}, vol.~35, no.~4, pp. 2693--2703, 2020.

\bibitem{cao2020deep}
J.~Cao, D.~Harrold, Z.~Fan, T.~Morstyn, D.~Healey, and K.~Li, ``Deep reinforcement learning-based energy storage arbitrage with accurate lithium-ion battery degradation model,'' \emph{IEEE Transactions on Smart Grid}, vol.~11, no.~5, pp. 4513--4521, 2020.

\bibitem{yousuf}
Y.~Baker, N.~Zheng, and B.~Xu, ``Transferable energy storage bidder,'' \emph{IEEE Transactions on Power Systems}, 2023.

\bibitem{li2015efficiency}
R.~Li, A.~J. Svoboda, and S.~S. Oren, ``Efficiency impact of convergence bidding in the california electricity market,'' \emph{Journal of Regulatory Economics}, vol.~48, pp. 245--284, 2015.

\bibitem{li2021machine}
Y.~Li, N.~Yu, and W.~Wang, ``Machine learning-driven virtual bidding with electricity market efficiency analysis,'' \emph{IEEE Transactions on Power Systems}, vol.~37, no.~1, pp. 354--364, 2021.

\bibitem{krishnamurthy2017energy}
D.~Krishnamurthy, C.~Uckun, Z.~Zhou, P.~R. Thimmapuram, and A.~Botterud, ``Energy storage arbitrage under day-ahead and real-time price uncertainty,'' \emph{IEEE Transactions on Power Systems}, vol.~33, no.~1, pp. 84--93, 2017.

\bibitem{Bitar}
E.~Y. Bitar, R.~Rajagopal, P.~P. Khargonekar, K.~Poolla, and P.~Varaiya, ``Bringing wind energy to market,'' \emph{IEEE Transactions on Power Systems}, vol.~27, no.~3, pp. 1225--1235, 2012.

\bibitem{xu2018optimal}
B.~Xu, Y.~Shi, D.~S. Kirschen, and B.~Zhang, ``Optimal battery participation in frequency regulation markets,'' \emph{IEEE Transactions on Power Systems}, vol.~33, no.~6, pp. 6715--6725, 2018.

\bibitem{Nie}
Y.~Nie, N.~H. Nguyen, P.~Sinthong, and J.~Kalagnanam, ``A time series is worth 64 words: Long-term forecasting with transformers,'' \emph{arXiv preprint arXiv:2211.14730}, 2022.

\bibitem{zhou2021informer}
H.~Zhou, S.~Zhang, J.~Peng, S.~Zhang, J.~Li, H.~Xiong, and W.~Zhang, ``Informer: Beyond efficient transformer for long sequence time-series forecasting,'' in \emph{Proceedings of the AAAI conference on artificial intelligence}, vol.~35, no.~12, 2021, pp. 11\,106--11\,115.

\bibitem{Vaswani}
A.~Vaswani, N.~Shazeer, N.~Parmar, J.~Uszkoreit, L.~Jones, A.~N. Gomez, {\L}.~Kaiser, and I.~Polosukhin, ``Attention is all you need,'' \emph{Advances in neural information processing systems}, vol.~30, 2017.

\bibitem{zeng}
A.~Zeng, M.~Chen, L.~Zhang, and Q.~Xu, ``Are transformers effective for time series forecasting?'' in \emph{Proceedings of the AAAI conference on artificial intelligence}, vol.~37, no.~9, 2023, pp. 11\,121--11\,128.

\bibitem{wu}
H.~Wu, J.~Xu, J.~Wang, and M.~Long, ``Autoformer: Decomposition transformers with auto-correlation for long-term series forecasting,'' \emph{Advances in Neural Information Processing Systems}, vol.~34, pp. 22\,419--22\,430, 2021.

\bibitem{hyndman}
R.~J. Hyndman and A.~B. Koehler, ``Another look at measures of forecast accuracy,'' \emph{International journal of forecasting}, vol.~22, no.~4, pp. 679--688, 2006.

\bibitem{zhao2007framework}
J.~H. Zhao, Z.~Y. Dong, X.~Li, and K.~P. Wong, ``A framework for electricity price spike analysis with advanced data mining methods,'' \emph{IEEE Transactions on Power Systems}, vol.~22, no.~1, pp. 376--385, 2007.

\bibitem{baltaoglu2018algorithmic}
S.~Baltaoglu, L.~Tong, and Q.~Zhao, ``Algorithmic bidding for virtual trading in electricity markets,'' \emph{IEEE Transactions on Power Systems}, vol.~34, no.~1, pp. 535--543, 2018.

\end{thebibliography}

% that's all folks
\end{document}